\documentclass[12pt]{article}

\usepackage{amsfonts}
\usepackage{amssymb}
\usepackage{amsthm}
\usepackage{epsf}
\usepackage{pstricks,pst-node,pst-text,pst-3d}
\usepackage{amsmath,stmaryrd}
\usepackage[mathscr]{eucal}
\usepackage{tikz}
\usetikzlibrary{patterns}
\usepackage[english]{babel}
\usepackage{subfig}
\usepackage{mathtools}
\usepackage{bm}
\usetikzlibrary{decorations.pathreplacing}
\usepackage{sectsty}
\usepackage{multirow}
\usepackage{url}
\usepackage{setspace}
\usepackage{enumitem}
\usepackage[utf8]{inputenc}

\sectionfont{\fontsize{14}{15}\selectfont}

\makeatletter
\DeclareRobustCommand{\rvdots}{%
  \vbox{
    \baselineskip4\p@\lineskiplimit\z@
    \kern-\p@
    \hbox{.}\hbox{.}\hbox{.}
  }}
\makeatother

\newtheorem{definition}{Definition}
\newtheorem{theorem}{Theorem}
\newtheorem{proposition}{Proposition}

\newtheorem{example}{Example}
\newtheorem{corollary}{Corollary}
\newtheorem{lemma}{Lemma}

\newcommand{\niton}{\not\owns}

\begin{document}

\title{\Large\bf On Admissibility in \\ Bipartite Incidence Graph Sampling}
\author{Pedro Garc\'ia-Segador \\ National Statistics Institute (Spain) \and Li-Chun Zhang\\University of Southampton, UK \\ Statistisk sentralbyrå, Norway}
\maketitle

\begin{abstract}
In bipartite incidence graph sampling, the target study units may be formed as connected population elements, which are distinct to the units of sampling and there may exist generally more than one way by which a given study unit can be observed via sampling units. This generalizes finite-population element or multistage sampling, where each element can only be sampled directly or via a single primary sampling unit. We study the admissibility of estimators in bipartite incidence graph sampling and identify other admissible estimators than the classic Horvitz-Thompson estimator. Our admissibility results encompass those for finite-population sampling. 
\end{abstract}


\section{Introduction}

In the development of finite-population sampling theory since Neyman (1934), the estimator of Horvitz and Thompson (1952) has a central position. Whilst uniformly minimum variance unbiased estimator (UMVUE) does not exist generally (e.g., Godambe, 1955; Hanurav, 1966; Basu, 1971), the Horvitz-Thompson estimator (HTE) is the only UMVUE in the class of linear unbiased estimators for a special type of sampling designs, called unicluster designs. However, unicluster designs (such as systematic sampling) do not satisfy all the needs of efficiency or practicability in applications, and it is not possible to estimate the sampling variance unbiasedly for these designs (Hanurav, 1966).    
 
Godambe (1960) shifts the attention to admissibility as an alternative criterion and proves that the HTE is admissible in the class of linear unbiased estimators. The linearity restriction is removed subsequently (Godambe and Joshi, 1965; Joshi, 1965, 1966), and the admissibility of HTE is extended to the entire class of unbiased estimators. Ramakrishnan (1973) gives an alternative, simpler proof of this result. 

In \emph{bipartite incidence graph sampling} (Zhang, 2021a, 2021b), shorthanded as BIGS, the sampling units are formally distinguished from the target study units and there can exist \emph{more than one way} by which a given study unit may be observed via sampling units. In other words, let a bipartite digraph have two sets of nodes representing  sampling and study units, respectively, such that directed edges only exist from the first set of nodes to the other, where each edge represents the incidence observation relationship between a sample unit and a study unit (regardless the actual operations involved or if the relationships may be unknown in advance). For example, sampling `influencers' on a social media platform via the edges from an initial sample of users is a case of BIGS, where an edge exists from any user to each one she `follows' and an influencer can be sampled via more than one user. Similarly, if one samples webpages by following the links from other webpages, where the study units are the webpages with links from others and the sampling units are all the webpages.

Such \emph{multiplicity of access} to a given study unit distinguishes BIGS from `standard' finite-population element or multistage sampling. For instance, in two-stage sampling (e.g., Cochran, 1977), each element can be sampled via one and only one primary sampling unit; similarly for each element called ultimate sampling unit in multistage sampling.

On the one hand, BIGS can extend the scope of estimation by allowing the study units to be formed as groups of population elements, such as social networks of individuals (e.g., Goodman, 1961; Frank, 1971, 1980) or a cluster of neighbouring habitats (Thompson, 1990). On the other hand, it enables a unified treatment (Zhang, 2021b) of so-called `non-standard' population sampling techniques as special cases of BIGS, such as multiplicity or indirect sampling (Birnbaum and Sirken, 1965; Lavallée, 2007), network sampling (Sirken, 2005) and adaptive cluster sampling (Thompson, 1990), as well as breadth-first or depth-first techniques that are explicitly devised for graph data, such as snowball sampling (Goodman, 1961; Frank and Snijders, 1994; Zhang and Oguz-Alper, 2020; Zhang and Patone, 2017) and random-walk type sampling (e.g., Thompson, 2006; Avrachenkov et al., 2010; Zhang, 2021, 2024). 

 A large family of unbiased \emph{incidence weighting estimators (IWEs)} has been proposed for BIGS (Zhang, 2021b; Patone and Zhang, 2022), which includes the HTE as a special case and many others. There is thus a need to study admissibility in BIGS, since there may be other admissible estimators than the HTE which may or may not be IWEs. 

The paper continues as follows: in Section 2, we introduce formally the notation, BIGS, IWE and other basic concepts. In Section 3, we prove the two main results of this work regarding the admissibility of unbiased estimators in BIGS, which encompass the existing results in finite-population sampling mentioned earlier. A summary is given in Section 4, together with some open problems for future research.

\section{Basic concepts}

Denote by $\mathcal{B}=(F,\Omega;H)$ a bipartite simple digraph, where $F$ is the non-empty node set of \emph{sampling units}, $\Omega$ the non-empty set of \emph{study units} and $H$ the set of \emph{edges} each of which points from a node in $F$ to a node in $\Omega$. Whenever possible, we will denote the elements of $F$ by $\{ i_1, i_2, ... \}$ and the elements of $\Omega$ by $\{ \kappa_1, \kappa_2, ...\}$. Let an initial sample $s_0$ be taken from $F$, $s_0 \subseteq F$. The nodes in $\Omega$ connected to those in $s_0$ are called \emph{successors} of $s_0$, denoted by $\alpha(s_0)$. For any $i\in F$, the set of successors of $i$ in $\Omega$ is denoted by $\alpha_i$. For a subset $\Lambda$ of $\Omega$, $\Lambda \subseteq \Omega$, the nodes in $F$ connected to $\Lambda$ are called the \emph{ancestors} of $\Lambda$, denoted by $\beta(\Lambda)$, and the set of ancestors of any $\kappa \in \Omega$ is denoted by $\beta_{\kappa}$. 

For BIGS given any initial sample $s_0$, the study units are given by the \emph{incident observation procedure}, denoted by $\Omega_s =\alpha(s_0)$, together with the corresponding edges emanating from $s_0$, $H_s=H\cap (s_0 \times \Omega)$. The \emph{sample graph} (from $\mathcal{B}$) is given as $\mathcal{B}_s=(s_0,\Omega_s; H_s)$. Sometimes, we will also use the notation $\Omega_s(s_0)$ to make explicit the choice of $s_0$. Note that any standard setting of finite-population element or multistage sampling can be represented as BIGS from a graph $\mathcal{B}$ where $|\beta_{\kappa}| \equiv 1$ for any element $\kappa$ in the population $\Omega$. 

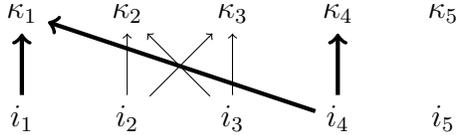
\begin{figure}[ht]
  \centering
\begin{tikzpicture}[scale=.7]

  \node (a1) at (0,0) {$i_1$};
  \node (a2) at (2,0) {$i_2$};
  \node (a3) at (4,0) {$i_3$};
  \node (a4) at (6,0) {$i_4$};
  \node (a5) at (8,0) {$i_5$};
  
  \node (b1) at (0,2) {$\kappa_1$};
  \node (b2) at (2,2) {$\kappa_2$};
  \node (b3) at (4,2) {$\kappa_3$};
  \node (b4) at (6,2) {$\kappa_4$};
  \node (b5) at (8,2) {$\kappa_5$};

  \draw [-to, ultra thick] (a1) -- (b1);
  \draw [-to] (a2) -- (b2);
  \draw [-to] (a2) -- (b3);
  \draw [-to] (a3) -- (b2);
  \draw [-to] (a3) -- (b3);
  \draw [-to, ultra thick] (a4) -- (b1); 
  \draw [-to, ultra thick] (a4) -- (b4); 

\end{tikzpicture}
\caption{Example of a bipartite graph}
\label{Fig_example1}
\end{figure}

Figure \ref{Fig_example1} provides an illustration. For instance, the successors of $i_2$ are $\alpha_{i_2}=\{ \kappa_2, \kappa_3 \}$, and the ancestors of $\kappa_1$ are $\beta_{\kappa_1}=\{ i_1, i_4 \}$. Note that there may exist unconnected elements in $F$ (such as $i_5$) or $\Omega$ (such as $\kappa_5$). If the initial sample is $s_0=\{ i_1, i_4 \}$ then $\Omega_s=\{ \kappa_1, \kappa_4\}$ and the edges $H_s$ of the sample graph $\mathcal{B}_s$ are marked in bold. 

Let $y_{\kappa}$ be an unknown constant associated with each study unit $\kappa\in \Omega$. Let $y(\Omega) =\{y_{\kappa} : \kappa \in \Omega \}$, or as a vector $y\in \mathbb{R}^{\vert \Omega \vert}$. BIGS from $\mathcal{B}$ allows us to observe $y(\Omega_s) = \{ y_{\kappa} : \kappa \in \Omega_s \}$ associated with the realised sample graph $\mathcal{B}_s$. Let the target of estimation be the total
\[
\theta = \sum_{\kappa \in \Omega}y_{\kappa}.
\]

For a given initial sampling design $p(s_0)$ on $F$, let $\mathcal{S}_0$ be the set of possible initial samples with $p(s_0)>0$ which is called the \emph{support} of $p$, such that $\sum_{s_0\in \mathcal{S}_0} p(s_0)=1$. Let the initial-sample inclusion probabilities be $\pi_i = \Pr(i \in s_0)$, $\pi_{ij} = \Pr(i \in s_0, j \in s_0)$ and so on, which we distinguish from the study-sample inclusion probabilities written as $\pi_{(\kappa)} = \Pr(\kappa \in \Omega_s)$ or $\pi_{(\kappa \ell)} = \Pr(\kappa \in \Omega_s, \ell\in \Omega_s)$ which are induced from $p(s_0)$ and the given graph $\mathcal{B}$. Due to the multiplicity of access to study units in BIGS, additional knowledge of $\mathcal{B}$ beyond $\mathcal{B}_s$ is needed to calculate $\pi_{(\kappa)}$ or $\pi_{(\kappa\ell)}$. We shall distinguish two levels of knowledge in this paper. 
\begin{itemize}[leftmargin=6mm]
\item \emph{Graph knowledge:} when the entire graph $\mathcal{B} = (F, \Omega; H)$ is known. 
\item \emph{Ancestry knowledge:} when only the ancestor sets $\mathfrak{AK}_s=\{ \beta_{\kappa} : \kappa \in \Omega_s \}$ are known for $\Omega_s$.
\end{itemize}

\begin{figure}[ht]
 \centering
\begin{tikzpicture}[scale=.7]
  \node (a1) at (0,0) {$i_1$};
  \node (a2) at (2,0) {$i_2$};
  \node (a3) at (4,0) {$i_3$};
  \node (a4) at (6,0) {$i_4$};
  \node (a5) at (8,0) {$i_5$};
  
  \node (b1) at (0,2) {$\kappa_1$};
  \node (b2) at (2,2) {$\kappa_2$};
  \node (b3) at (4,2) {$\kappa_3$};
  \node (b4) at (6,2) {$\kappa_4$};
  \node (b5) at (8,2) {$\kappa_5$};
  
  \node (c1) at (0,3) {$y_{\kappa_1}$};
  \node (c2) at (2,3) {$y_{\kappa_2}$};
  \node (c3) at (4,3) {$y_{\kappa_3}$};

  \draw [-to, ultra thick] (a1) -- (b1);
  \draw [-to, ultra thick] (a2) -- (b1);
  \draw [-to, ultra thick] (a2) -- (b2);
  \draw [-to, ultra thick] (a2) -- (b3);
  \draw [-to, ultra thick] (a3) -- (b3);
  \draw [-to, ultra thick] (a3) -- (b4);
  \draw [-to, ultra thick] (a4) -- (b4); 
  \draw [-to, ultra thick] (a4) -- (b2); 
  \draw [-to, ultra thick] (a5) -- (b5); 


  \node (a1) at (12,0) {$i_1$};
  \node (a2) at (14,0) {$i_2$};
  \node (a3) at (16,0) {$i_3$};
  \node (a4) at (18,0) {$i_4$};
  \node (a5) at (20,0) {$i_5$};
  
  \node (b1) at (12,2) {$\kappa_1$};
  \node (b2) at (14,2) {$\kappa_2$};
  \node (b3) at (16,2) {$\kappa_3$};
  \node (b4) at (18,2) {$\kappa_4$};
  \node (b5) at (20,2) {$\kappa_5$};
  
  \node (c1) at (12,3) {$y_{\kappa_1}$};
  \node (c2) at (14,3) {$y_{\kappa_2}$};
  \node (c3) at (16,3) {$y_{\kappa_3}$};
  
  \draw [-to, ultra thick] (a1) -- (b1);
  \draw [-to, ultra thick] (a2) -- (b1);
  \draw [-to, ultra thick] (a2) -- (b2);
  \draw [-to, ultra thick] (a2) -- (b3);
  \draw [-to, ultra thick] (a3) -- (b3);
  \draw [-to] (a3) -- (b4);
  \draw [-to] (a4) -- (b4); 
  \draw [-to, ultra thick] (a4) -- (b2); 
  \draw [-to] (a5) -- (b5); 
\end{tikzpicture}
 \caption{Graph knowledge (left), ancestry knowledge (right) given $s_0 = \{ i_1, i_2\}$.}
 \label{Fig_example_knowledge}
\end{figure}
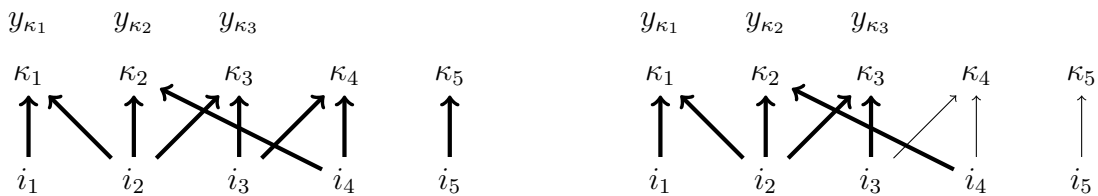

Figure \ref{Fig_example_knowledge} illustrates the differences between the two levels of knowledge, given the initial sample $s_0=\{ i_1, i_2 \}$, where the edges highlighted in bold represent the known part of the graph. Note that at least the ancestry knowledge is needed for the calculation of $\pi_{(\kappa)}$ and, consequently, the HTE. Meanwhile, even in situations where the graph knowledge may be available in principle, such as the linked webpages earlier, the graph may be too large or non-static over time such that sampling can be more resource-effective to obtain an estimate of $\theta$ than to count $\theta$ directly by attempting to process the whole graph. 

Now we can give a formal definition of estimator in BIGS.

\begin{definition}\label{def:estimator}
An \emph{estimator} $e(s_0; y , \mathcal{B})$ is a real-valued function on $\mathcal{S}_0 \times \mathbb{R}^{\vert \Omega \vert} \times \mathfrak{B}$, where $\mathfrak{B}$ denotes all possible $\mathcal{B}$ given $(F, \Omega)$, which depends on $y$ only through $\{ y_{\kappa} : \kappa\in \alpha(s_0)\}$ and depends on $\mathcal{B}$ only through $\mathcal{B}_s$ and the given level of knowledge. We may use the simplified notation $e$, $e(s_0)$ or $e(s_0; \mathcal{B})$ to refer to an estimator when the context is clear.
\end{definition}

Notice that, given graph knowledge, the condition that an estimator $e$ only depends on the values of $y$ through the study sample translates to
\[
e(s_0;y,\mathcal{B})=e(s_0;y',\mathcal{B}),\quad \forall y(\Omega_s) = y'(\Omega_s) 
\]
i.e. for every $y'$ that coincides with $y$ in the realised $\Omega_s$. Whereas, given ancestry knowledge, the conditions that an estimator $e$ must meet would translate to
\[
e(s_0;y,\mathcal{B})=e(s_0;y',\mathcal{B}'),\quad \forall y(\Omega_s) = y'(\Omega_s),~\mathcal{B}_s = \mathcal{B}'_s,~ \mathfrak{AK}_s = \mathfrak{AK}_s'
\]
i.e. for every $y'$ that coincides with $y$ in $\Omega_s$ \emph{and} every $\mathcal{B}'$ that coincides with $\mathcal{B}$ in terms of the realised same sample graph and the corresponding ancestor sets of $\Omega_s$.

Next, throughout this work, we shall be working under the condition that every $y_{\kappa}$ in $\Omega$ has a positive probability to be sampled, which is formally phrased as a property of $\mathcal{B}$. 

\begin{definition}
A graph $\mathcal{B}=(F,\Omega;H)$ is \emph{covered} by a design $p(s_0)$ on $F$ if 
\[
\pi_{(\kappa)}= \Pr(\kappa \in \Omega_s) >0, \ \forall \kappa \in \Omega.
\]
Let us denote by $\mathfrak{B}(p)$ the collection of graphs $\mathcal{B}$ covered by a given $p(s_0).$
\end{definition}

Given $p(s_0)$, we denote by $D$ the class of all unbiased estimators of $\theta$ based on BIGS, where an estimator $e$ is unbiased if, over hypothetically repeated sampling, 
\[
\text{E}\{ e(s_0; y, \mathcal{B}) \} = \theta 
\]
for any given $\mathcal{B} \in \mathfrak{B}(p)$ and associated $y\in \mathbb{R}^{\vert \Omega \vert}$. A general approach to unbiased estimators is provided by the family of IWEs (Zhang, 2021b), given as
\[
\hat{\theta}_{IWE} =\sum_{(i\kappa) \in H_s} W_{i \kappa}\dfrac{y_{\kappa}}{\pi_i}
=\sum_{\kappa \in \Omega_s} \Big( \sum_{i \in \beta_{\kappa}\cap s_0} \dfrac{W_{i\kappa}}{\pi_i} \Big)  y_{\kappa}
=\sum_{i \in s_0} \dfrac{1}{\pi_i} \Big( \sum_{\kappa \in \alpha_{i}} W_{i\kappa} y_{\kappa} \Big),
\]
where the weights $W_{i \kappa}$ can vary with $\mathcal{B}_s$ generally. Given any $\mathcal{B} \in \mathfrak{B}(p)$, the IWE $\hat{\theta}_{IWE}$ is unbiased iff for each $\kappa \in \Omega$, we have 
\begin{equation}\label{eq:insesga}
\sum_{i \in \beta_{\kappa}} \text{E}(W_{i \kappa} \vert i\in s_0)=1.
\end{equation}
(Zhang, 2021b, Theorem 2.1). In particular, if the weight assigned to every $(i\kappa)\in H$ is a constant, regardless of $\mathcal{B}_s$ and denote by $w_{i\kappa}$ for distinction, the IWE is unbiased iff for each $\kappa \in \Omega$, we have 
\[
\sum_{i \in \beta_{\kappa}} w_{i \kappa}=1.
\]

Various unbiased estimators in non-standard finite-population sampling applications can be studied as special cases of the IWE; see Zhang (2021b) and Patone and Zhang (2022). The following examples of unbiased IWE will be used in the discussions in this paper. Firstly, the HTE defined on $\Omega_s$ can be written as an IWE, i.e.
\begin{equation} \label{eq:HTE}
\hat{\theta}_{HTE}=\sum_{\kappa \in \Omega_s} \dfrac{y_{\kappa}}{\pi_{(\kappa)}} = \hat{\theta}_{IWE} \qquad\text{given}\quad
 \sum_{i \in s_0 \cap \beta_{\kappa}} \dfrac{W_{i \kappa}}{\pi_i}=\dfrac{1}{\pi_{(\kappa)}}
\end{equation}
as a condition on the weights $W_{i \kappa}$. Moreover, let $s_{\kappa} =s_0 \cap \beta_{\kappa}$ and $\phi_{s_{\kappa}} =\Pr(s_0\cap \beta_{\kappa}=s_{\kappa})$, a \emph{HT-type estimator} is given by the weights $W_{i \kappa}$ satisfying 
\[
\eta_{s_{\kappa}}=\pi_{(\kappa)} \sum_{i \in s_{\kappa}} \dfrac{W_{i \kappa}}{\pi_i} \qquad\text{and}\qquad
\sum_{s_{\kappa}} \phi_{s_{\kappa}}\eta_{s_{\kappa}} = \pi_{(\kappa)}
\]
which includes the HTE \eqref{eq:HTE} as the special case of $\eta_{s_{\kappa}} \equiv 1$. Secondly, the so-called \emph{Hansen-Hurwitz type (HH-type) estimator} is an IWE that uses constant weights $w_{i \kappa}$. The earliest example is the \emph{multiplicity estimator} (Birnbaum and Sirken, 1965) using equal weights
\begin{equation} \label{equal-w}
w_{i \kappa}=\vert \beta_{\kappa} \vert ^{-1}.
\end{equation}
Note that all these IWEs are feasible given the ancestry knowledge of $\Omega_s$. 

Let us now give a definition of admissibility in classes of unbiased estimators in BIGS. 

\begin{definition}
Let $p(s_0)$ be a design on $F$ and $\mathcal{D}$ a class of unbiased estimators, $\mathcal{D}\subseteq D$. An estimator $e_0 \in \mathcal{D}$ is said to be \emph{admissible} in $\mathcal{D}$, if there is no other $e \in \mathcal{D}$ satisfying
\[
\begin{cases} V(e) \leq V(e_0), & \text{for \emph{all}}~ y(\Omega), \ \mathcal{B}\in \mathfrak{B}(p) \\ V(e) < V(e_0), & \text{for \emph{some}}~ y(\Omega), \ \mathcal{B}\in \mathfrak{B}(p) \end{cases}
\] 
In case an estimator is not admissible, it is said to be \emph{inadmissible}.
\end{definition}

We note immediately that if an estimator $e_0$ is admissible in $\mathcal{D}$, then there cannot exist another estimator in $\mathcal{D}$ which has the same variance given any $(\mathcal{B}, y)$. The result is given below and its proof in Appendix \ref{sec:proof} (as are all the other proofs in this paper).

\begin{proposition}\label{prop:strongly_admi}
Let $\mathcal{D}\subseteq D$ be a class of unbiased estimators such that if $e_0,e_0+d\in \mathcal{D},$ then $e_0+ a d \in \mathcal{D}, \ \forall a \in \mathbb{R}$. If $e_0$ is admissible in the class $\mathcal{D}$, then it is the only estimator $e\in \mathcal{D}$ such that
\[
\text{V} \left( e \right) \leq \text{V} \left( e_0 \right), \ \forall \mathcal{B} \in \mathfrak{B(p)},\, y(\Omega).
\]
\end{proposition}

Note that the hypothesis of $\mathcal{D}$ holds when $\mathcal{D}=D$, since $e_0,e_0+d\in D \Leftrightarrow \text{E}(d)=0,$ thus $\text{E}(ad) =a\text{E}(d) = 0$ and $e_0+ a d\in D$. All the classes of estimators we will work with in this paper will meet this hypothesis of $\mathcal{D}$.

However, the whole class $D$ is still too broad to establish admissibility \emph{generally}, although we already know that the HTE is admissible in $D$. Let us therefore introduce an additional property of the estimators that is necessary to our development later.

\begin{figure}[ht]
\centering
\begin{tikzpicture}[scale=.7]  
  \node (c0) at (8,-1) {$\mathcal{B}$};
  \node (c1) at (6,0) {$i_1$};
  \node (c2) at (8,0) {$i_2$};
  \node (c3) at (10,0) {$i_3$};
  
  \node (d1) at (6,2) {$\kappa_1$};
  \node (d2) at (8,2) {$\kappa_2$};
  \node (d1m) at (6,3) {$0$};
  \node (d2m) at (8,3) {$y_{\kappa_2}$};

  \draw [-to] (c1) -- (d1);
  \draw [-to] (c1) -- (d2);
  \draw [-to] (c2) -- (d2);
  \draw [-to] (c3) -- (d2);
  
  \node (e0) at (14,-1) {$\mathcal{B}^{(\kappa_1)}$};
  \node (e1) at (12,0) {$i_1$};
  \node (e2) at (14,0) {$i_2$};
  \node (e3) at (16,0) {$i_3$};
  
  \node (f1) at (14,2) {$\kappa_2$};
  \node (f2m) at (14,3) {$y_{\kappa_2}$};

  \draw [-to] (e1) -- (f1);
  \draw [-to] (e2) -- (f1);
  \draw [-to] (e3) -- (f1);
\end{tikzpicture}
\caption{A zero-invariant estimator is the same given $\mathcal{B}$ or $\mathcal{B}^{(\kappa_1)}$.}
 \label{Fig_zero}
\end{figure}
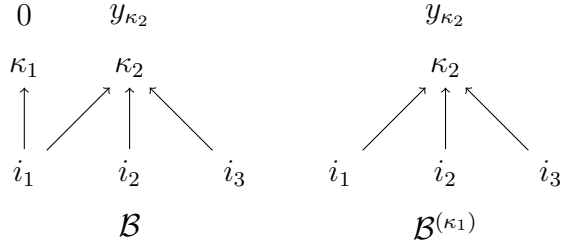

Let $\mathcal{B}=(F,\Omega; H)$ be a graph and $\Lambda \subseteq \Omega$, we denote by $\mathcal{B}^{(\Lambda)}$ the graph resulting from removing from $\mathcal{B}$ the nodes $\Lambda$ and any edge in $H$ that ends in $\Lambda$. When the set $\Lambda$ consists of a single element $\kappa$, we will use the notation $\mathcal{B}^{(\kappa)}$. Now, in the case $y_{\kappa} \equiv 0$ for any $\kappa\in \Lambda$, the target total $\theta$ is the same in $\mathcal{B}$ and $\mathcal{B}^{(\Lambda)}$, as illustrated in Figure \ref{Fig_zero}, and it seems intuitive that the differences between $\mathcal{B}$ and $\mathcal{B}^{(\Lambda)}$ should not affect the expectation and variance of any potentially admissible estimator of $\theta$. This is formalized in the definition below, and the potential loss of efficiency of non-zero-invariant estimators will be demonstrated later.

\begin{definition}
An estimator $e$ is \emph{zero-invariant} if for any $s_0\in \mathcal{S}_0$, we have
\[
e(s_0; y, \mathcal{B}) = e(s_0; y', \mathcal{B}^{(\Lambda)})
\]
for any $y(\Omega)$ and $y'(\Omega \setminus \Lambda)$ satisfying $y_{\kappa}=0$ if $\kappa \in \Lambda$ and $y_{\kappa} =y_{\kappa}'$ if $\kappa \in \Omega \setminus \Lambda$. In the special case of empty $\Omega_s$, $\Omega_s=\emptyset$, we will let the estimator take the value zero. Denote by $D^*$ the class of \emph{unbiased zero-invariant estimators}.
\end{definition}

Note that the HTE is zero-invariant, since $y_{\kappa}/\pi_{(\kappa)} = 0$ if $y_{\kappa} = 0$ regardless $\pi_{(\kappa)}$. In this work we study admissibility for classes of unbiased zero-invariant estimators.

\section{Admissibility in classes of unbiased estimators}

Below we consider admissibility given either graph knowledge or ancestry knowledge.

\subsection{Given graph knowledge}

Given the graph knowledge, one can determine the study-sample $\Omega_s(s_0)$ following any initial sample $s_0$, so as to apply Rao-Blackwellization (Rao, 1945; Blackwell, 1947) and the concept of sufficiency (Godambe and Joshi, 1965). Indeed, we shall prove below that any sufficient estimator is admissible in the class $D^*$ of unbiased zero-invariant estimators, and any such estimator can be obtained by a modified Rao-Blackwellization procedure. 

Recall that in finite population sampling by $p(s)$, yielding sample $s$ from population $U$, the minimal sufficient statistic is the unordered sample $\{(i, y_i) : i \in s\}$, with all $\{ y_i : i\in U\}$ as the unknown parameters. In BIGS that yields the study sample $\Omega_s(s_0)$, denote the minimal sufficient statistic by
\[
D_{s_0}:=\{ (\kappa, y_{\kappa}) : \kappa \in \Omega_s(s_0) \}.
\]
We shall characterise an estimator as \emph{sufficient} by the following result without proof.  

\begin{proposition}\label{prop:sufi}
In BIGS, an estimator $e$ is \emph{sufficient} for $y$ given $p(s_0)$ and any $\mathcal{B}$, if
\[
e(s_0) = e(s'_0)
\]
for $s_0 \neq s'_0$ as long as $\Omega_s(s_0)=\Omega_s(s'_0)$.
\end{proposition}

Notice that any sufficient estimator can as well be written as $e(\Omega_s)$ with $\Omega_s$ as the \emph{distinct} images of all the possible initial samples, including $\Omega_s = \emptyset$, where the sampling probability of $\Omega_s$ follows from $p(s_0)$, denoted by
\[
P(\Omega_s) = \Pr\{ \alpha(s_0) = \Omega_s\} = \sum_{s_0 : \alpha(s_0) = \Omega_s} p(s_0)
\] 
Clearly, given any estimator $e$, one can apply Rao-Blackwellization to obtain a sufficient estimator, which is given as
\[
e_{RB} = \text{E}(e \mid D_{s_0}) = \sum_{s_0' : \alpha(s_0') = \Omega_s(s_0)} \frac{p(s_0') e(s_0')}{P(\Omega_s(s_0))}
\]
and we would obtain $e_{RB} = e$ if the estimator $e$ is sufficient, such as the HTE. 

Although Rao-Blackwellization does not yield UMVUE because $D_{s_0}$ is not complete in BIGS (nor in sampling from $\Omega$ directly), we shall prove below that it is closely related to admissibility. However, it is necessary to modify the procedure since $e_{RB}$ is not zero-invariant generally, which can have undesirable consequences, as the example below illustrates.

\begin{example}\label{example2}
Consider again the graphs in Figure \ref{Fig_zero}. Let there be two possible initial samples $s_0^* = \{i_1, i_2 \}$ and $s_0^{**} = \{i_2, i_3 \}$, where $p(s_0^*)=1/3$ and $p(s_0^{**})=2/3$. 
We have $\pi_1=1/3$, $\pi_2 =1$ and $\pi_3=2/3$. Let $e$ be the multiplicity estimator \eqref{equal-w}, where
\[
\begin{cases} e(s_0^*; y', \mathcal{B}^{(\kappa_1)}) = (\pi_1^{-1} + \pi_2^{-1}) \frac{y_{\kappa_2}}{3} = \frac{4}{3} y_{\kappa_2} \\
e(s_0^{**}; y', \mathcal{B}^{(\kappa_1)}) = (\pi_2^{-1} + \pi_3^{-1}) \frac{y_{\kappa_2}}{3} = \frac{5}{6} y_{\kappa_2} \end{cases} \text{or}\quad
\begin{cases} e(s_0^*; y, \mathcal{B}) = 3y_{\kappa_1} + \frac{4}{3} y_{\kappa_2} \\
e(s_0^{**}; y, \mathcal{B}) = \frac{5}{6} y_{\kappa_2} \end{cases}
\]
given $\mathcal{B}^{(\kappa_1)}$ or $\mathcal{B}$. The corresponding $e_{RB}$ is 
\[
e_{RB}(s_0^*; y', \mathcal{B}^{(\kappa_1)}) = e_{RB}(s_0^{**}; y', \mathcal{B}^{(\kappa_1)}) = y_{\kappa_2}
\]
given $\mathcal{B}^{(\kappa_1)}$; whereas 
\[
e_{RB}(s_0^*; y, \mathcal{B}) = e(s_0^*; y, \mathcal{B}) \not\equiv e(s_0^{**}; y, \mathcal{B}) = e_{RB}(s_0^{**}; y, \mathcal{B}) 
\]
because $\alpha(s_0^*) \neq \alpha(s_0^{**})$ given $\mathcal{B}$. Thus, given $y_{\kappa_1} =0$ and $\mathcal{B}$ instead of $\mathcal{B}^{(\kappa_1)}$, the variation of $e_{RB}$ over $s_0$ causes extra variance unnecessarily. Note that being zero-invariant, the HTE is always equal to $y_{\kappa_2}$ here given $\mathcal{B}^{(\kappa_1)}$ or $\mathcal{B}$.
\end{example}

Let use therefore modify Rao-Blackwellization so that it becomes zero-invariant. 

\begin{definition}
The \emph{zero-invariant Rao-Blackwell (ZRB)} estimator derived from any given unbiased estimator $e$ is
\[
e_{ZRB}(s_0) \coloneqq \text{ZRB}(e) = \sum_{s_0' \in [s_0]} \frac{p(s_0') e(s_0')}{P(\Omega_s(s_0) \cap \Omega_{(0)})}
\]
where $\Omega_{(0)} = \{ \kappa\in \Omega : y_{\kappa} \neq 0\}$, and $\Omega_s(s_0) \cap \Omega_{(0)}$ contains only the non-zero study units in the realised sample $\Omega_s(s_0)$, and $[s_0] = \{ s_0' : \alpha(s_0') \cap \Omega_{(0)} = \Omega_s(s_0) \cap \Omega_{(0)} \}$.
\end{definition}

By definition, $e_{ZRB}$ is zero-invariant if $e$ is zero-invariant, and it is unbiased if $e$ unbiased. Moreover, $e_{ZRB}$ has a smaller mean squared error if $e$ is not constant over $[s_0]$, as illustrated with $\mathcal{B}$ in Example \ref{example2} above. The operation, denoted by $\text{ZRB}(e)$, achieves thus the intended effect of Rao-Blackwellization in BIGS. The results are formalized below, where the theorem is the main result of admissibility given graph knowledge.

\begin{lemma}\label{lema:lagrangiano}
Let $p$ be a design and any $\mathcal{S}_0^* \subseteq \mathcal{S}_0$. The minimum of $\sum_{s_0 \in \mathcal{S}_0^*} p(s_0) e^2(s_0)$ subjected to the constraint $\sum_{s_0 \in \mathcal{S}_0^*} p(s_0) e(s_0) = Y^*$ is only attained at
\[
e(s_0) = \dfrac{Y^*}{\sum_{s_0 \in \mathcal{S}_0^*} p(s_0)}.
\]
\end{lemma}

\begin{proposition} \label{prop:RB}
Let $e$ be an unbiased zero-invariant estimator in BIGS. The corresponding $e_{RB}$ is inadmissible in $D$, since 
\[
\text{V}(e_{ZRB}) \leq \text{V}(e_{RB}) \leq \text{V}(e),
\]
where the first inequality is strict for some graph $\mathcal{B}$ and vector $y$.
\end{proposition}

\begin{theorem}\label{th:main_sufficient}
Assume graph knowledge. Let $p(s_0)$ be a design, where $\pi_i >0$ for any $i\in F$ and $e_0$ an estimator in the class $D^*$. Then the following statements are equivalent:
\begin{itemize}
\item[$i)$] $e_0$ is sufficient;
\item[$ii)$] $e_0$ is admissible in the class $D^*$;
\item[$iii)$] $e_0$ is given by $\text{ZRB}(e)$ for some estimator $e \in D^*$.
\end{itemize}
\end{theorem}

We note that the graph knowledge is needed for $\text{ZRB}(e)$ but not for the proof of (i) and (ii) in Theorem \ref{th:main_sufficient}. For instance, the sufficient, zero-invariant and admissible HTE requires only the ancestry knowledge. Moreover, now that BIGS encompasses finite-population sampling, Theorem \ref{th:main_sufficient} applies directly to all the non-standard population sampling problems, which clarifies the connection of admissibility to sufficiency and Rao-Blackwellization. 

For example, Thompson (1990) proposes the \emph{modified HTE} for adaptive cluster sampling and applies Rao-Blackwellization to it, where the units with $y$-values below the fixed adaptive threshold are used in the estimator \emph{only when} they are \emph{directly selected} in the initial sample, not otherwise when they are observed via the units with $y$-values above the threshold. The modified HTE is thus zero-invariant but not sufficient, which is why it can be modified by Rao-Blackwellization. Theorem \ref{th:main_sufficient} shows that this indeed yields an admissible estimator in the class $D^*$ for adaptive cluster sampling; indeed, applying ZRB to any unbiased zero-invariant estimator would yield an admissible estimator for adaptive cluster sampling.

Returning to any unbiased IWE estimator that is zero-invariant due to its linearity in $y$, applying ZRB to it would yield an admissible estimator in $D^*$. The admissible HTE is a special case that is zero-invariant and sufficient, such that it is unchanged by ZRB.  

\begin{corollary}
Assume graph knowledge. Let $p(s_0)$ be a design, where $\pi_i >0$ for any $i\in F$. The estimator $e_{ZRB}$ derived from some unbiased IWE $e$ is admissible in the class $D^*$. 
\end{corollary}

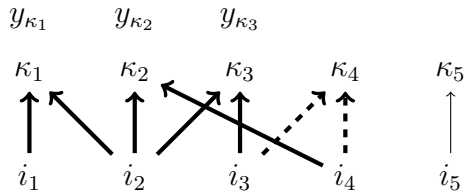
\begin{figure}[ht]
\centering
\begin{tikzpicture}[scale=.7]
  \node (a1) at (0,0) {$i_1$};
  \node (a2) at (2,0) {$i_2$};
  \node (a3) at (4,0) {$i_3$};
  \node (a4) at (6,0) {$i_4$};
  \node (a5) at (8,0) {$i_5$};
  
  \node (b1) at (0,2) {$\kappa_1$};
  \node (b2) at (2,2) {$\kappa_2$};
  \node (b3) at (4,2) {$\kappa_3$};
  \node (b4) at (6,2) {$\kappa_4$};
  \node (b5) at (8,2) {$\kappa_5$};
  
  \node (c1) at (0,3) {$y_{\kappa_1}$};
  \node (c2) at (2,3) {$y_{\kappa_2}$};
  \node (c3) at (4,3) {$y_{\kappa_3}$};

  \draw [-to, ultra thick] (a1) -- (b1);
  \draw [-to, ultra thick] (a2) -- (b1);
  \draw [-to, ultra thick] (a2) -- (b2);
  \draw [-to, ultra thick] (a2) -- (b3);
  \draw [-to, ultra thick] (a3) -- (b3);
  \draw [-to, ultra thick, dashed] (a3) -- (b4);
  \draw [-to, ultra thick, dashed] (a4) -- (b4); 
  \draw [-to, ultra thick] (a4) -- (b2); 
  \draw [-to] (a5) -- (b5);  
\end{tikzpicture}
 \caption{Ancestry knowledge (bold) and successor-ancestry knowledge (bold and dashed), given $s_0 = \{ i_1, i_2\}$. Graph knowledge comprises all edges (see also Figure \ref{Fig_example_knowledge}) regardless $s_0$.}
 \label{Fig_example_sak}
\end{figure}

A final remark on the assumed graph knowledge is worthwhile. While the condition is easy to state and sufficient for Theorem \ref{th:main_sufficient}, what is actually needed is the part of the graph that enables one to apply ZRB to the \emph{realised} study sample $\Omega_s$. It is then not necessary to know the whole graph, but it suffices to know all the successors of any node in the ancestor set $\mathfrak{AK}_s$. One may refer to this level of knowledge as the \emph{successor-ancestry knowledge}, which can replace the assumption of graph knowledge in order to apply Theorem \ref{th:main_sufficient}. Figure \ref{Fig_example_sak} illustrates the different levels of knowledge when the initial sample is $s_0=\{ i_1, i_2 \}$.

\subsection{Given ancestry knowledge}

When we only have the ancestry knowledge, the HTE estimator is still admissible because it is sufficient and zero-invariant. However, for other estimators that are not sufficient, we can no longer invoke ZRB as by Theorem \ref{th:main_sufficient}. We therefore need to identify some additional properties of the estimators which are shared by the HTE but not limited to it. 

The first additional property we require for the main result to be developed concerns analytic functions. A mathematical function is analytic at a point if its local behavior can be expressed as an infinite sum of powers (i.e. its Taylor series) converging in a region around that point. An important property of nonzero analytic functions $f(x)$ is that their zeros are isolated, which means there is a neighborhood around each zero $x_0$, $f(x_0) =0$, where the function does not vanish except at $x_0$. 

We shall call an estimator \emph{analytic} if it is an analytic function of $y(\Omega_s)$. For example, every IWE estimator is linear and thus analytic, whereas $\text{ZRB}(e)$ is not generally continuous at zero and therefore is not analytic. Below we will denote by $D^{**}$ \emph{the class of unbiased, zero-invariant and analytic estimators} and study admissibility in $D^{**}$.

Next, we will call an estimator \emph{elemental} if it takes the same expression as the HTE for a special type of graphs, called elemental graphs, as formalised below.

\begin{definition}
In BIGS, an estimator $e$ is \emph{elemental} if $e(s_0,\mathcal{B})= y_{\kappa}/\pi_i$ whenever $i \in s_0$, given any \emph{elemental} graph $\mathcal{B}=(F,\Omega;H)$ with $\Omega=\{ \kappa \} $ and $H=\{ (i,\kappa)\}$.
\end{definition} 

\begin{figure}[ht]
\centering
\begin{tikzpicture}[scale=.7]

\node (a1) at (0,0) {$i_1$};
\node (a2) at (2,0) {$i_2$};
\node (a3) at (4,0) {$i_3$};
\node (b1) at (4,2) {$\kappa$};

\draw [-to] (a3) -- (b1);
\end{tikzpicture}
\caption{Example of elemental graph.}
\label{Fig_regu}
\end{figure}
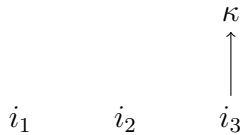

Figure \ref{Fig_regu} shows an example of elemental graph. We would have $e(s_0; y,\mathcal{B}) = y_{\kappa}/\pi_3$ given any sample $s_0$ containing $i_3$. Note that if we assume an estimator $e$ is linear, elemental and zero-invariant, then for any graph where $\beta_{\kappa}=1,$ $\forall \kappa \in \Omega,$ this estimator will match the HTE in BIGS, which will also coincide with the HTE by sampling from $\Omega$ directly. Being elemental is thus a property of the HTE which may be shared by other estimators such as the multiplicity estimator.

Next, we will invoke some linear algebra concepts that are relevant to elemental estimators and the proof of the main result later.

\begin{definition} Let the \emph{sample-space matrix} $\mathbb{S}^p$ have rows corresponding to the sampling units $F$ and columns to all the possible samples in the support of $p(s_0)$, $s_0 \in \mathcal{S}_0.$ The element for row $i$ and column $s_0$ is $\mathbb{S}^p_{i,s_0} =1$ if $i\in s_0$ and $0$ otherwise. Likewise, we define the matrix $\tilde{\mathbb{S}}^p$ whose element is $\tilde{\mathbb{S}}^p_{i,s_0} =p(s_0)$ if $i\in s_0$ and $0$ otherwise. 
\end{definition}

\begin{definition} An estimator $e$ is \emph{sample-space spanned} if there exist functions $a_i$ that do not depend on (i.e. vary with) the sample $s_0$ such that:
\[
e(s_0; y, \mathcal{B}) = \sum_{i\in s_0} a_i(y,\mathcal{B}), \quad\forall s_0 \in \mathcal{S}_0.
\]
\end{definition}

\begin{proposition}\label{prop: spanned-elemental}
Let $e \in D^*.$ If $e$ is sample-space spanned then $e$ is elemental.
\end{proposition}

Now, given a linear map defined by matrix $A$, the \emph{kernel} of $A$, $\text{Ker}(A)$, is the vector space formed by all the vectors $v$, such that $A\cdot v=0$. The \emph{row space}, $\text{Row}(A)$, is the linear space spanned by the row vectors of $A$. By definition, an estimator $e$ is sample-space spanned if and only if $e\in \text{Row}(\mathbb{S}^p)$ given any $y$ and $\mathcal{B}$. Note that $\text{Row}(A)$ is orthogonal to $\text{Ker}(A)$, denoted by $\text{Ker}(A)^{\perp}=\text{Row}(A)$. 
Moreover, let $v \cdot_p w = \sum_i p_i v_i w_i$ be the $p$-weighted scalar product of vectors $v$ and $w$, and $W^{\perp_p}$ the orthogonal complement of a vector space $W$, such that $v\cdot_p w=0$ for any $v\in W^{\perp_p}$, $w\in W$. We can now state the following result.

\begin{lemma}\label{lem:tech_alg}
Let $p(s_0)$ be a design with associated matrices $\mathbb{S}^p$ and $\tilde{\mathbb{S}}^p$. We have
\[
\text{Ker}(\tilde{\mathbb{S}}^p)^{\perp_p} = \text{Row}(\mathbb{S}^p).
\]
\end{lemma}

It follows from Lemma \ref{lem:tech_alg} that 
\begin{equation}\label{eq:direct_sum}
\text{Row}(\mathbb{S}^p) \cap \text{Ker}(\tilde{\mathbb{S}}^p)=\bm{0} \quad\text{and}\quad 
\mathbb{R}^{\vert \mathcal{S}_0 \vert}=\text{Row}(\mathbb{S}^p) \oplus \text{Ker}(\tilde{\mathbb{S}}^p)
\end{equation}
where $\bm{0}$ is the null-vector and $\oplus$ denotes direct sum. However, the decomposition \eqref{eq:direct_sum} of vector space $\mathbb{R}^{\vert \mathcal{S}_0 \vert}$ does not imply that every estimator in BIGS can be given as the sum of an estimator that is sample-space spanned and another that lies in the kernel. To make this distinction we introduce the next definition.

\begin{definition}
Let $p(s_0)$ be a design and $V \subseteq \mathbb{R}^{\vert \mathcal{S}_0 \vert}$ a vector space. We define the associated \emph{estimator space}, $\mathfrak{E}(V)$, as the set of all estimators $e$ such that $e(s_0;y,\mathcal{B}) \in V, \ \forall (y, \mathcal{B})$. Denote by $\mathfrak{E}_c(V)$ the subset of estimators in $ \mathfrak{E}(V)$ whose expectation is $c$ over sampling.
\end{definition}

If $V \cap W=\bm{0}$, then $\mathfrak{E}(V) \oplus \mathfrak{E}(W)$ contains the estimators that can be uniquely expressed as the sum of an estimator in $V$ and another in $W$. In particular, $\mathfrak{E}_{\theta}\big( \text{Row}(\mathbb{S}^p) \big)$ is just the set of unbiased sample-space spanned estimators of $\theta$, and $\mathfrak{E}_0\big(\text{Ker}(\tilde{\mathbb{S}}^p) \big)$ is the set of estimators $d$ in the specified kernel space satisfying $\text{E}(d) = 0$. Their relevance to delineating admissibility is clarified in the following proposition and used in the proof of the theorem next, which states that every estimator in $\mathfrak{E}_{\theta}\big( \text{Row}(\mathbb{S}^p) \big)$ is admissible in $D^{**}$.

\begin{proposition}\label{prop: ortho-inadmi}
Let $p(s_0)$ be a design and $e =e_0+d$ an estimator, where $e_0\in \mathfrak{E}_{\theta}\big( \text{Row}(\mathbb{S}^p) \big)$ and $d\in \mathfrak{E}_0\big(\text{Ker}(\tilde{\mathbb{S}}^p) \big)$. Then $\text{V}(e_0) \leq \text{V}(e)$ for every $(y, \mathcal{B})$, with strict inequality whenever $d\not\equiv 0$ for given $(y, \mathcal{B})$, i.e., if $d$ is not constant 0. In particular, if $d\not\equiv 0,$ $e$ is inadmissible.
\end{proposition}

\begin{theorem}\label{th: main_non_sufficient}
Assume ancestry knowledge in BIGS and let $p$ be a design for which $\pi_i>0,$ $\forall i \in F$. If $e_0$ in $D^{**}$ is a sample-space spanned estimator, i.e. $e_0 \in \mathfrak{E}_{\theta}\big( \text{Row}(\mathbb{S}^p) \big)$, then $e_0$ is admissible in the class $D^{**}.$ 
\end{theorem}

Note that sample-space spanned estimators can be seen as a special case of estimators in $\mathfrak{E}_{\theta}\big( \text{Row}(\mathbb{S}^p) \big) \oplus \mathfrak{E}_0\big(\text{Ker}(\tilde{\mathbb{S}}^p) \big)$ with the trivial kernel element $d\equiv 0$. Hence, we can summarize the results obtained so far for the estimators in $D^{**}$ as follows:
\begin{itemize}[leftmargin=6mm]
\item any estimator $e_0$ in $\mathfrak{E}_{\theta}\big( \text{Row}(\mathbb{S}^p) \big)$ is admissible in $D^{**}$, whereas
any estimator $e= e_0 + d$ with $e_0\in \mathfrak{E}_{\theta}\big( \text{Row}(\mathbb{S}^p) \big)$ and non-trivial $d \in \mathfrak{E}_0\big(\text{Ker}(\tilde{\mathbb{S}}^p) \big)$  is inadmissible in $D^{**}$;
\item the admissibility of estimator $e$ is unknown generally if $e\notin \mathfrak{E}_{\theta}\big( \text{Row}(\mathbb{S}^p) \big) \oplus \mathfrak{E}_0\big(\text{Ker}(\tilde{\mathbb{S}}^p) \big)$.
\end{itemize}
In particular, any IWE with constant weights $w_{i\kappa}$ is a linear (hence, analytic) sample-space spanned estimator, which can be given as $\hat{\theta}_{IWE}=\sum_{i \in s_0} a_i(y,\mathcal{B})$ where
\[
a_i(y,\mathcal{B}) = \sum_{\kappa \in \alpha_i} w_{i\kappa} y_{\kappa}/\pi_i. 
\]
Any such unbiased IWE is admissible in $D^{**}$ as summarised below without proof. 

\begin{corollary}
Assume ancestry knowledge and let $p(s_0)$ be a design for which $\pi_i>0,$ $\forall i \in F$. An unbiased IWE with constant weights $w_{i\kappa}$ is admissible in the class $D^{**}.$ 
\end{corollary}

Let us close this section with two further remarks. First, although the admissibility of estimators that lie outside of $\mathfrak{E}_{\theta}\big( \text{Row}(\mathbb{S}^p) \big) \oplus \mathfrak{E}_0\big(\text{Ker}(\tilde{\mathbb{S}}^p) \big)$ is yet unknown generally, it is still possible to obtain specific admissibility results by making restrictions on the designs $p(s_0)$, similarly to the aforementioned result that the HTE is the only UMVUE given unicluster designs in finite-population sampling. Below we give such a result.

\begin{definition} A design $p(s_0)$ over $F$ is said to be \emph{full rank} if $\mathbb{S}^p$ is a full rank matrix, that is, if $\text{rank}(\mathbb{S}^p) = |\mathcal{S}_0|$. 
\end{definition}

It is not difficult to verify that, for designs with fixed sample size $|s_0|$, all minimum support designs (e.g. Tillé, 2006) are full rank. Moreover, systematic sampling designs are full rank, but the simple random sampling design is not full rank. It now follows from Theorem \ref{th: main_non_sufficient} that every estimator in $D^{**}$ is admissible given full-rank designs.

\begin{corollary} \label{cor_fullrank}
Assume ancestry knowledge and let $p(s_0)$ be a full-rank design for which $\pi_i >0$, $\forall i \in F$. Every estimator in $D^{**}$ is admissible in the class $D^{**}.$
\end{corollary}

Second, the construction of Proposition \ref{prop: ortho-inadmi} to inadmissible estimators can be generalized to orthogonal pairs of estimator spaces other than $\mathfrak{E}_{\theta}\big( \text{Row}(\mathbb{S}^p) \big) \oplus \mathfrak{E}_0\big(\text{Ker}(\tilde{\mathbb{S}}^p) \big)$, denoted by $\mathfrak{E}_{\theta}(V) \oplus \mathfrak{E}_0(W)$ given
$W = V^{\perp_p}$, where $\mathfrak{E}_{\theta}(V)$ is a space of unbiased estimators of $\theta$ and $\mathfrak{E}_0(W)$ its orthogonal space of mean-zero estimators. Below we give such an example. 

\begin{example}\label{ex: example_final}
Consider simple random sampling design $p(s_0) \equiv 1/6$ with $|s_0| =2$ given the graph in Figure \ref{Fig_example1}, where we remove the disconnected elements $i_5$ and $\kappa_5.$ Let an IWE $e_0$ with non-constant weights $W_{i\kappa}$ be given as follows. Denote by $12<13<14<23<24<34$ the lexicographical order among the 6 samples in $\mathcal{S}_0$, based on which we define the weights $W_{i\kappa}$ and $e_0$ to be 
\[
W_{i\kappa}(s_0) = \dfrac{w_{i\kappa} \pi_i I_{i}(s_0)}{p(s_0)} \quad\text{and}\quad
e_0 =\sum_{\kappa \in \Omega_s} \left( \sum_{i \in \beta_{\kappa}\cap s_0} \dfrac{w_{i\kappa} I_{i}(s_0)}{p(s_0)} \right)  y_{\kappa}
\]
where $w_{i\kappa}$ can be any constant weights satisfying $\sum_{i\in \beta_{\kappa}} w_{i\kappa} =1$, $\forall \kappa$, and $I_{i}(s_0)$ is an indicator that $s_0$ is the first sample in the lexicographical order where $s_0 \ni i$. Since $I_i(s_0) =1$ only for one particular $s_0$ and $I_i(s_0) =0$ otherwise, each $y_{\kappa}$ is used for estimation given only one particular $s_0$ as well; hence the expression of $e_0$. The estimator $e_0$ is unbiased because
\[
\sum_{i \in \beta_{\kappa}} \text{E}(W_{i \kappa} \vert i \in s_0) = \sum_{i \in \beta_{\kappa}} \sum_{s_0 \ni i} \dfrac{p(s_0) W_{i \kappa}}{\pi_i}=\sum_{i \in \beta_{\kappa}} w_{i\kappa}  \sum_{s_0 \ni i}  I_{i}(s_0)=\sum_{i \in \beta_{\kappa}} w_{i\kappa}=1 .
\]
Note that $e_0(s_0) =0$ when $s_0 = \{ i_2, i_3\}$, $\{ i_2, i_4\}$ or $\{ i_3, i_4\}$, where $I_{i}(s_0)=0$ for any $i\in s_0$. In other words, the estimator $e_0$ resides in the vector space 
\[
V= \{ (a,b,c,0,0,0): a,b,c \in \mathbb{R} \}.
\]

Next, to construct a mean-zero estimator $d$ that belongs to the orthogonal vector space, let us repeat the construction of $e_0$ with different ordered samples. Let $e_1$ be given as $e_0$ except that $I_i(s_0)$ is based on the reverse lexicographical order $34<24<23<14<13<12$; and let $e_2$ be based on the order $24<34<23<14<13<12$. As both $e_1$ and $e_2$ are unbiased, their difference $d=e_2 - e_1$ has zero mean, which resides in the vector space
\[
W = \{ (0,0,0,0,f,g): f,g \in \mathbb{R} \}.
\]
Table \ref{tab:example_final} lists the estimates $e_0$, $e_1$, $e_2$ and $d$ given each of the 6 samples $s_0$, where we use the multiplicity weights $w_{i\kappa}=1/\vert \beta_{\kappa} \vert$ for this illustration.

\begin{table}[ht]
\centering
\caption{Estimators $e_0$, $e_1$, $e_2$ and $d$ given $w_{i\kappa} = 1/\vert \beta_{\kappa}\vert$.}
\small
\setlength{\tabcolsep}{4pt}
\begin{tabular}{c|c|c|c|c} \hline
 $s_0$ & $e_0(s_0)$ & $e_1(s_0)$ & $e_2(s_0)$ & $d(s_0)$ \\ \hline
$\{ i_1, i_2\}$ & $3(y_{\kappa_1}+y_{\kappa_2}+y_{\kappa_3})$ & $0$ & $0$ & $0$ \\ \hline
$\{ i_1, i_3\}$ & $3(y_{\kappa_2}+y_{\kappa_3})$ & $0$ & $0$ & $0$ \\ \hline
$\{ i_1, i_4\}$ & $3(y_{\kappa_1}+2y_{\kappa_4})$ & $3y_{\kappa_1}$ & $3y_{\kappa_1}$ & $0$ \\ \hline
$\{ i_2, i_3\}$ & $0$ & $0$ & $0$ & $0$ \\ \hline
$\{ i_2, i_4\}$ & $0$ & $3(y_{\kappa_2}+y_{\kappa_3})$ & $3(y_{\kappa_1}+y_{\kappa_2}+y_{\kappa_3}+2y_{\kappa_4})$ & $3(y_{\kappa_1}+2y_{\kappa_4})$ \\ \hline
$\{ i_3, i_4\}$ & $0$ & $3(y_{\kappa_1}+y_{\kappa_2}+y_{\kappa_3}+2y_{\kappa_4})$ & $3(y_{\kappa_2}+y_{\kappa_3})$ & $-3(y_{\kappa_1}+2y_{\kappa_4})$ \\ \hline
\end{tabular} \label{tab:example_final}
\end{table}

Clearly, $V$ and $W$ are orthogonal, but $\mathfrak{E}_{\theta}(V)\oplus \mathfrak{E}_0(W) \neq \mathfrak{E}_{\theta}\big( \text{Row}(\mathbb{S}^p) \big) \oplus \mathfrak{E}_0\big(\text{Ker}(\tilde{\mathbb{S}}^p) \big)$, since $e_0$ is not a sample-space spanned estimator. One can now follow the proof of Proposition \ref{prop: ortho-inadmi} in Appendix \ref{sec:proof} to show that $\text{V}(e_0) < \text{V}(e_0 +d)$; hence, $e=e_0+d$ is inadmissible.
\end{example}

\section{Concluding remarks}

In this work, we have analyzed the admissibility of many unbiased estimators in BIGS, depending on the available knowledge of the graph, and identified many other admissible estimators in BIGS than the HTE. 

In the case of graph knowledge, the investigation has naturally led us to the concepts of sufficiency, Rao-Blackwellization and its zero-invariant modification ZRB. In Theorem \ref{th:main_sufficient}, we have proved that for unbiased zero-invariant estimators, admissibility is equivalent to sufficiency and the application of ZRB. Thus, given graph knowledge (or successor-ancestry knowledge which enables the application of ZRB for the realised sample graph), it will always be better to perform ZRB of any unbiased estimator to ensure its admissibility. 

In the case of ancestry knowledge, we have established in Theorem \ref{th: main_non_sufficient} that all sample-space spanned estimators are admissible in the class of unbiased, zero-invariant and analytic estimators, as well as devised a constructive approach to identify inadmissible estimators therein by means of orthogonal vector spaces (Proposition \ref{prop: ortho-inadmi} and Example \ref{ex: example_final}). 

\begin{table}[h!]
\centering
\caption{Admissibility of some IWE estimators in BIGS}
\begin{tabular}{c|c|c|c|c} \hline
Knowledge & Estimator class & HTE & HT-type & HH-type \\ \hline
Graph & $D^*$ & Admissible & Inadmissible & Inadmissible \\ \hline
Ancestry & $D^{**}$ & Admissible & Unknown generally & Admissible \\ \hline
\end{tabular} \label{tab:estimadores}
\end{table}

The results regarding the examples of IWE introduced at the beginning are summarized in Table \ref{tab:estimadores}. As one can see, with graph knowledge, the HTE is admissible because it is sufficient, but the others are inadmissible because they are not sufficient. However, provided ancestry knowledge, the HH-type estimators are admissible in the class $D^{**}$ and cannot be dominated by the HTE; whereas the HT-type estimators with non-constant weights $W_{i \kappa}$ are not covered by Theorem \ref{th: main_non_sufficient}, so their admissibility is still unknown generally, although there are special cases such as Corollary \ref{cor_fullrank} if one only considers restricted designs.

Finally, since all the sampling-unbiased estimators known in the literature are analytic, restricting to the class $D^{**}$ is not a limitation that matters much for the practice of sampling theory, compared to studying the classes $D^*$ or $D$ more broadly. An obvious topic for future research would be to further delineate admissibility for the estimators in $D^{**}$, possibly by other means than pairwise orthogonal vector spaces devised in this paper.

\appendix 
\section{Proofs} \label{sec:proof}

Proof of Proposition \ref{prop:strongly_admi}.
\begin{proof}
Let us assume that $e$ meets the condition of the statement. Since no strict inequality can occur because $e_0$ is admissible, the variance must always be equal $\text{V} (e) = \text{V} (e_0)$, that is, if $e=e_0+d,$ then
\[
\text{V}(e) = \text{V}(e_0) + \text{V}(d) + 2\text{Cov}(e_0,d) \Leftrightarrow \text{V}(d) = -  2\text{Cov}(e_0,d).
\]
Observe that $\text{Cov}(e_0,d)\leq 0$. If for some graph $\mathcal{B}$ and vector $y$, $e \neq e_0 \Rightarrow d\neq 0$ and $\text{Cov}(e_0,d)<0$, then, since $e_\alpha=e_0+\alpha d \in  \mathcal{D}$ for any $\alpha \in (0,1)$,  we get
\[
\text{V} (e_\alpha) = \text{V} (e_0) + \alpha^2 \text{V} (d) + 2\alpha\text{Cov} (e_0,d) =  \text{V} (e_0) + 2\alpha(1-\alpha)\text{Cov} (e_0,d) \leq \text{V} (e_0),
\]
and $\text{V} (e_\alpha) < \text{V} (e_0)$ for the cases where  $\text{Cov} (e_0,d)<0 .$ This is a contradiction with the fact that $e_0$ is admissible, such that $e=e_0$ is the only possibility.
\end{proof}

\noindent
Proof of Lemma \ref{lema:lagrangiano}.
\begin{proof}
By the Lagrangian of the optimization problem, we have
\begin{gather*}
L(e,\lambda)=\sum_{s_0 \in \mathcal{S}_0^*} p(s_0) e^2(s_0) - \lambda \left( \sum_{s_0 \in \mathcal{S}_0^*} p(s_0) e(s_0) - Y^* \right), \\
\dfrac{\partial L}{\partial e(s_0)}= 2 p(s_0) e(s_0) - \lambda p(s_0) =0 \quad\Rightarrow\quad e(s_0) = \dfrac{\lambda}{2 }, \ \ \forall s_0 \in \mathcal{S}_0^*.
\end{gather*}
Using the constraint we obtain
\[
\lambda=\dfrac{2 Y ^*}{\sum_{s_0 \in \mathcal{S}_0^*} p(s_0)} \quad\Rightarrow\quad
e(s_0) = \dfrac{\lambda}{2 }= \dfrac{Y^*}{\sum_{s_0 \in \mathcal{S}_0^*} p(s_0)}.
\]
Since the Hessian is positive definite, the solution corresponds to a local minimum. As we are minimizing a convex function over a convex set, the minimum is also global. Moreover, since the function is strictly convex, the minimum is unique and the lemma holds.
\end{proof}

\noindent
Proof of Proposition \ref{prop:RB}.
\begin{proof}
The second inequality is simply the Rao-Blackwell theorem. Let us focus on the first inequality. If $y(\Omega)$ are all non-zero, then $e_{ZRB}=e_{RB}$. Otherwise, let $y(\Lambda)$ consist of all the zero values, where $\Lambda \subset \Omega$. Between $\mathcal{B}^{(\Lambda)}$ and $\mathcal{B}$ let us define
\[
\mathfrak{O}(\mathcal{B}, \Lambda) = \{ \Omega_s(s_0, \mathcal{B}^{(\Lambda)}) : \exists s_0, s'_0 \ \text{s.t. } \alpha(s_0; \mathcal{B})\neq \alpha(s'_0; \mathcal{B}) \text{ and } \alpha(s_0; \mathcal{B}^{(\Lambda)}) = \alpha(s'_0; \mathcal{B}^{(\Lambda)})  \}
\]
and work by induction on the cardinality $K = |\mathfrak{O}(\mathcal{B}, \Lambda)|$. 

First, if $K =1$, let the only element of $\mathfrak{O}(\mathcal{B}, \Lambda)$ be $\Omega_s(s'_0; \mathcal{B}^{(\Lambda)}) = \Omega_s(s'_0; \mathcal{B}) \cap \Omega_{(0)}$ , and $[s'_0] = \{ s_0 :  \alpha(s_0; \mathcal{B}^{(\Lambda)})  = \Omega_s(s'_0; \mathcal{B}^{(\Lambda)})\}$. Since both $e_{RB}$ and $e_{ZRB}$ are unbiased, and equal to each other whenever $s_0\notin [s'_0],$ the problem reduces to comparing the second moments for the elements in $[s'_0].$ Additionally, note that the estimators $e, e_{RB}$ and $e_{ZRB}$ satisfy:
 
\[
\sum_{s_0 \in [s'_0]} p(s_0) e(s_0) = Y^* := \theta- \sum_{s_0 \notin [s'_0]} p(s_0)e(s_0).
\]
Applying Lemma \ref{lema:lagrangiano} with $\mathcal{S}_0^* = [s_0']$ yields $V(e_{ZRB}) \leq V(e_{RB})$ where the inequality is strict if $|[s'_0]| > 1$. 

Next, suppose $K >1$. Define estimators $e_{ZRB,1}, e_{ZRB,1}, ..., e_{ZRB,K} = e_{ZRB}$, where $e_{ZRB,1}$ is given by applying ZRB only to the first element of $\mathfrak{O}(\mathcal{B}, \Lambda)$, and $e_{ZRB,2}$ by applying ZRB only to the first two elements of $\mathfrak{O}(\mathcal{B}, \Lambda)$, and so on. Applying Lemma \ref{lema:lagrangiano} to compare $e_{ZRB,k}$ and $e_{ZRB,k-1}$, for $k= K, ..., 2$, as well as $e_{ZBR,1}$ and $e_{RB}$, we obtain
\[
\text{V}(e_{ZRB}) = \text{V}(e_{ZRB,K}) \leq \text{V}(e_{ZBR,K-1}) \leq \cdots \leq \text{V}(e_{ZBR,1}) \leq \text{V}(e_{RB}).
\]
This completes the proof. 
\end{proof}

\noindent
Proof of Theorem \ref{th:main_sufficient}.

\begin{proof}
$i) \Rightarrow ii)$ We shall use a double induction on $N=\vert F \vert$ and $K=\vert \Omega \vert.$ For $N=1$ the only possible design is $s_0=\{ 1 \}$ and $p(s_0)=1$. Since $e_0$ is unbiased, i.e. $p(s_0)e_0(s_0)=\theta$, the only possible estimator is $e_0(s_0)=\theta$ which is therefore trivially admissible. 

We now assume that the result is true for $N.$ To show that the theorem holds for $N+1,$ we will show that for any given design $p(s_0)$ given $|F| =N+1$, where $\pi_i>0$ for any $i \in F$, for any estimator $e(s_0; y,\mathcal{B}) \in D^*$ such that
\begin{equation}\label{eq:hipo}
\text{V}\left( e(s_0; y,\mathcal{B}) \right) \leq \text{V}\left( e_0(s_0; y,\mathcal{B}) \right), \ \ \ \forall \mathcal{B} \in \mathfrak{B}(p), \ y(\Omega) ,
\end{equation}
it holds that $e(s_0; y,\mathcal{B})= e_0(s_0; y,\mathcal{B})$ for every $s_0 \in \mathcal{S}_0$.

We will prove the above property by induction on $K=\vert \Omega \vert.$ If $K=1,$ the graphs associated with this case are in the form shown in Figure \ref{Fig_main_1}. There could be elements in $F$ that do not connect to $\kappa$, and if there is a sample formed exclusively by these elements, the value of the estimator, due to the zero-invariance property, would be zero.

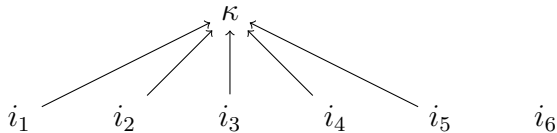
\begin{figure}[ht]
\centering
\begin{tikzpicture}[scale=.7]   
  \node (e1) at (0,0) {$i_1$};
  \node (e2) at (2,0) {$i_2$};
  \node (e3) at (4,0) {$i_3$};
  \node (e4) at (6,0) {$i_4$};
  \node (e5) at (8,0) {$i_5$};
  \node (e6) at (10,0) {$i_6$};
  
  \node (f1) at (4,2) {$\kappa$};

  \draw [-to] (e1) -- (f1);
  \draw [-to] (e2) -- (f1);
  \draw [-to] (e3) -- (f1);
  \draw [-to] (e4) -- (f1);
  \draw [-to] (e5) -- (f1);
\end{tikzpicture}
 \caption{Illustration of a graph in the case $K=1$.}
 \label{Fig_main_1}
\end{figure}

Note that either $\alpha(s_0)=\emptyset$ and the estimator takes the value zero, or $\alpha(s_0)=\{ \kappa \}.$ Since $e_0$ is sufficient, for any pair of samples $s_0, s'_0$ such that $\alpha(s_0)=\alpha(s_0')=\{ \kappa \}$ the estimator takes the same value, i.e. $e_0(s_0; y,\mathcal{B})=e_0(s'_0,y,\mathcal{B}).$ As the estimator is unbiased, we have:
\begin{gather*}
\sum_{\alpha(s_0)=\{ \kappa \}} p(s_0) e_0(s_0; y,\mathcal{B}) = \theta \quad\Leftrightarrow\quad e_0(s_0; y,\mathcal{B}) \sum_{\alpha(s_0)=\{ \kappa \}} p(s_0) = \theta \\
\Leftrightarrow\quad e_0(s_0; y,\mathcal{B}) = \dfrac{\theta}{\sum_{\alpha(s_0)=\{ \kappa \}} p(s_0)},~ \forall s_0 \text{ s.t. } \alpha(s_0)=\{ \kappa \}.
\end{gather*}
If we apply Lemma \ref{lema:lagrangiano}, considering that $Y^*=\theta$ and $\mathcal{S}^*_0$ is the set of samples that satisfy $\alpha(s_0)=\{ \kappa \}$, we obtain that the estimator above corresponds to the one that minimizes the variance.  As $e$ meets inequality (\ref{eq:hipo}), then $e=e_0$ and the result holds for $K=1.$ 

Let us assume that the result is true for $K$ and prove it for $K+1.$ Let $\mathcal{B} \in \mathfrak{B}(p)$ such that $\vert F \vert = N+1$ and $\vert \Omega \vert = K+1$ and $\kappa' \in \Omega$. Let $y_{\kappa '}=0.$ Now consider the graph $\mathcal{B}^{(\kappa')}$ where we remove $\kappa '$ and any edge that ends at $\kappa'$. Since $e$ and $e_0$ are zero-invariant we have
\[
e(s_0,\mathcal{B})=e(s_0,\mathcal{B}^{(\kappa')})=e_0(s_0,\mathcal{B}^{(\kappa')})=e_0(s_0,\mathcal{B}),
\]
where the first and third equalities are due to zero-invariance and the second one due to the induction hypothesis since $\mathcal{B}^{(\kappa')}$ has $K$ study units. 

In this way, we have shown that for any value of $y$ that vanishes at some element $\kappa' \in \Omega$ and for any graph $\mathcal{B}$, $e(s_0; y,\mathcal{B})=e_0(s_0; y,\mathcal{B})$ for any $s_0 \in \mathcal{S}_0$. Therefore, considering any sample $s_0$ such that $\alpha(s_0) \neq \Omega $ and a general vector $y$, we can choose some $\kappa ' \in \Omega \setminus \alpha(s_0)$ and, by the same argument, show that $e(s_0; y,\mathcal{B})=e_0(s_0; y,\mathcal{B})$ for any graph.

Let us see what happens if the sample $s_0$ satisfies $\alpha(s_0) = \Omega $. In this case, as $e$ meets $e_0$ for the rest of samples $s_0'$ with $\alpha(s_0') \neq \Omega ,$ from unbiasedness and inequality (\ref{eq:hipo}) we get
\begin{gather}
\sum_{\alpha(s_0)=\Omega} p(s_0) e(s_0; y,\mathcal{B}) = Y^* = \theta - \sum_{\alpha(s_0) \neq \Omega} p(s_0) e_0(s_0; y,\mathcal{B}), \label{eq:fin1} \\
\sum_{\alpha(s_0)=\Omega} p(s_0) e^2(s_0; y,\mathcal{B}) \leq \sum_{\alpha(s_0)=\Omega} p(s_0) e^2_0(s_0; y,\mathcal{B}). \label{eq:fin2}
\end{gather}
Since $e_0$ is sufficient and $\alpha(s_0)=\Omega$  we get from (\ref{eq:fin1}):
\[
\sum_{\alpha(s_0)=\Omega} p(s_0) e_0(s_0; y,\mathcal{B}) = Y^* \quad\Rightarrow\quad e_0(s_0; y,\mathcal{B})= \dfrac{Y^*}{\sum_{\alpha(s_0)=\Omega} p(s_0)}.
\]
If we apply Lemma \ref{lema:lagrangiano} for minimizing the left-hand side of (\ref{eq:fin2}) subject to (\ref{eq:fin1}), considering $\mathcal{S}^*_0$ as the set of samples that satisfy $\alpha(s_0)=\Omega$, we obtain the above estimator as the one that minimizes the variance and then by (\ref{eq:fin2}), $e(s_0; y,\mathcal{B})=e_0(s_0; y,\mathcal{B})$ and the result holds.

$ii) \Rightarrow iii)$  Note that if $e_0$ is not $\text{ZRB}(e)$ of any estimator $e \in D^*$, then $e_0  \neq \text{ZRB}(e_0)$ and by Proposition \ref{prop:RB}, $ \text{ZRB}(e_0)$ improves upon $e_0$ for some graph. Thus, $e_0$ is not admissible.

Another way to prove this is by applying the Rao-Blackwell Theorem, which shows that $\text{V}(e_{ZRB})\leq \text{V}(e_0).$ If $e_0$ were admissible, we would have $e_0=ZRB(e_0)$ by Proposition \ref{prop:strongly_admi}.

$iii) \Rightarrow i)$ If $e_0 = \text{ZRB}(e)$ for some $e \in D^*$, then $e_0$ is a sufficient estimator in $D^*$.
\end{proof}

\noindent
Proof of Proposition \ref{prop: spanned-elemental}.
\begin{proof}
Let $B$ be an elemental graph with the only edge $H=\{(j, \kappa ) \}$. As $e$ is zero-invariant, we have $a_i =0$ if $i \neq j$, and $e(s_0; y, \mathcal{B}) = a_j(y_{\kappa}, \mathcal{B})$ for every $s_0$ containing $j$. It follows from unbiasedness that $a_j(y_{\kappa}, \mathcal{B}) = y_{\kappa}/\pi_i$, i.e., $e$ is elemental.
\end{proof}

\noindent
Proof of Lemma \ref{lem:tech_alg}.
\begin{proof}
Between weighted and unweighted scalar products, we have
\[
v\in \text{Ker}(\tilde{\mathbb{S}}^p)^{\perp_p} \iff \sum_{s_0 \in \mathcal{S}_0} p(s_0)v(s_0)d(s_0)=0, \ \forall d \in \text{Ker}(\tilde{\mathbb{S}}^p) 
\iff pv\in \text{Ker}(\tilde{\mathbb{S}}^p)^{\perp}. 
\] 
Since $\text{Ker}(\tilde{\mathbb{S}}^p)^{\perp} =\text{Row}(\tilde{\mathbb{S}}^p)$, we have $pv \in \text{Row}(\tilde{\mathbb{S}}^p) \iff v \in \text{Row}(\mathbb{S}^p)$.
\end{proof}

\noindent 
Proof of Proposition \ref{prop: ortho-inadmi}.
\begin{proof}
We have $\text{V}(e)=\text{V}(e_0) + \text{V}(d) + 2\text{Cov}(e_0,d)$. Since $\text{E}(e)=\theta$ and $\text{E}(d)=0$, we obtain 
\[
\text{Cov}(e_0,d)=\sum_{s_0}p(s_0)\{ e_0(s_0) - \theta \} d(s_0) = \sum_{s_0} p(s_0)e_0(s_0)d(s_0) = e_0 \cdot_p d = 0
\]
due to the orthogonality established by Lemma \ref{lem:tech_alg}, such that $\text{V}(e_0) \leq \text{V}(e)$ for every $(y, \mathcal{B})$, and $\text{V}(e_0) < \text{V}(e)$ whenever $\text{V}(d)>0$ because $d\not \equiv 0$ for given $(y, \mathcal{B})$.
\end{proof}

\noindent 
Proof of Theorem \ref{th: main_non_sufficient}.
\begin{proof}
Suppose there exists some estimator $e \in D^{**}$ such that 
\[
\text{V}\left( e(s_0; y,\mathcal{B}) \right) \leq \text{V}\left( e_0(s_0; y,\mathcal{B}) \right), \ \ \ \forall \mathcal{B} \in \mathfrak{B}_{N}, \ y \in \mathbb{R}^{\vert \Omega(\mathcal{B}) \vert }. 
\]
We shall prove that $e(s_0; y,\mathcal{B})=e_0(s_0; y,\mathcal{B})$, $\forall s_0 \in \mathcal{S}_0$, $y \in \mathbb{R}^{\vert \Omega(\mathcal{B}) \vert }$ and $\mathcal{B}\in \mathfrak{B}(p)$. 

To do so, we adopt the following approach. First, suppose that $e$ differs to $e_0$ given some graph $\mathcal{B}$. Note that, since they are both unbiased, we can express $e_0$ and $e$ given $\mathcal{B}$ as
\[
e_0(s_0; y, \mathcal{B}) = \theta + h_0(s_0; y, \mathcal{B}) \quad\text{and}\quad e(s_0; y, \mathcal{B}) = \theta + h(s_0; y, \mathcal{B})
\]
where $\sum_{s_0 \in \mathcal{S}_0} p(s_0) h_0(s_0) =\sum_{s_0 \in \mathcal{S}_0} p(s_0) h(s_0)=0$. 

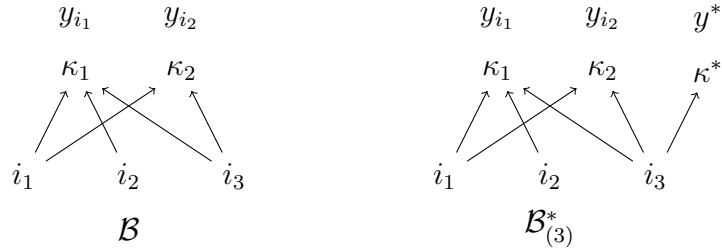
\begin{figure}[ht]
\centering
\begin{tikzpicture}[scale=.7]

  \node (a0) at (2,-1) {$\mathcal{B}$};
  \node (a1) at (0,0) {$i_1$};
  \node (a2) at (2,0) {$i_2$};
  \node (a3) at (4,0) {$i_3$};
  
  \node (b1) at (1,2) {$\kappa_1$};
  \node (b1m) at (1,3) {$y_{i_1}$};
  
  \node (b2) at (3,2) {$\kappa_2$};
  \node (b2m) at (3,3) {$y_{i_2}$};

  \draw [-to] (a1) -- (b1);
  \draw [-to] (a2) -- (b1);
  \draw [-to] (a3) -- (b1);
  
  \draw [-to] (a1) -- (b2);
  \draw [-to] (a3) -- (b2);
    
  \node (c0) at (10,-1) {$\mathcal{B}^*_{(3)}$};
  \node (c1) at (8,0) {$i_1$};
  \node (c2) at (10,0) {$i_2$};
  \node (c3) at (12,0) {$i_3$};
  
  \node (d1) at (9,2) {$\kappa_1$};
  \node (d1m) at (9,3) {$y_{i_1}$};
  
  \node (d2) at (11,2) {$\kappa_2$};
  \node (d2m) at (11,3) {$y_{i_2}$};  
  
  \node (d3) at (13,2) {$\kappa^*$};
  \node (d3m) at (13,3) {$y^*$}; 
  
  \draw [-to] (c1) -- (d1);
  \draw [-to] (c2) -- (d1);
  \draw [-to] (c3) -- (d1);
  \draw [-to] (c1) -- (d2);
  \draw [-to] (c3) -- (d2);
  
  \draw [-to] (c3) -- (d3);
\end{tikzpicture}
\caption{Illustration of $\mathcal{B}^*_{(i_3)}$ constructed from $\mathcal{B}$.}
\label{Fig_lema_non1}
\end{figure}

Next, let us construct a new graph $\mathcal{B}^*_{(j)}$ from $\mathcal{B}$ by adding a new element $\kappa^*$ to $\Omega$ and the edge $(j, \kappa^*)$ to $H$, with the value $y^*$ associated with $\kappa^*$, as illustrated in Figure \ref{Fig_lema_non1}. We shall prove that there exists some vector $\overline{y}$ associated with some $\mathcal{B}^*_{(j)}$, consisting of $y(\Omega)$ and $y^*$, where the variance of $e_0$ is less than that of $e$, contradicting the hypothesis that $\text{V}(e_0) \geq \text{V}(e)$ and resulting in the only possibility that $e$ and $e_0$ could not differ given $\mathcal{B}$. 

To argue for the contradiction, we shall distinguish whether $e$ is elemental or not, noting that $e_0$ is elemental by Proposition \ref{prop: spanned-elemental}.

\paragraph{Case 1:} Suppose $e$ is not elemental. Let us take $y=\bm{0}$ and compare the variance of $e_0$ and $e.$ As $e$ is not elemental there exists some $j^*\in F$ such that for the elemental graph $\mathcal{B}^*_{(j^*)}$ there is some sample $s^*_0$ such that $e(s^*_0; y^*, \mathcal{B}^*_{(j^*)})\neq y^*/\pi_{j^*}$. Now that both $e_0$ and $e$ are zero-invariant, we have $e_0(s_0) = e(s_0)=0$ for all $s_0 \niton j^*$. Moreover, since $e_0$ and $e$ are unbiased estimators, both satisfy

\[
\sum_{s_0 \in \mathcal{S}_0} p(s_0) e(s_0; \bar{y}, \mathcal{B}^*_{(j^*)}) = \sum_{s_0 \ni j^*} p(s_0) e(s_0; \bar{y}, \mathcal{B}^*_{(j^*)})= y^*.
\]

It follows that had $e(s_0)$ been the same for every sample $s_0 \ni j^*$, we could only have $e(s_0; y^*, \mathcal{B}^*_{(j^*)}) = y^*/\pi_{j^*}$ whenever $s_0 \ni j^*$. Therefore, $e(s_0)$ cannot be constant over $s_0$ if it is not elemental. Thus, by Lemma \ref{lema:lagrangiano}, the variance of $e_0$ achieves the minimum since $e_0$ is elemental, such that $\text{V}(e) > \text{V}(e_0)$, which yields the contradiction.

\paragraph{Case 2:} Suppose $e$ is elemental, i.e. $e(s_0; \bar{y}, \mathcal{B}^*_{(j)}) = y^*/\pi_j$ whenever $s_0\ni j$ given any choice of $j\in F$ in the case $y= \bm{0}$. As the estimators $e$ and $e_0$ are elemental and analytic, they are, in particular, continuous, and the following holds:
\[
\lim_{y \to \bm{0}} e_0(s_0; \bar{y},\mathcal{B}^*_{(j)})= \lim_{y \to \bm{0}} e(s_0; \bar{y},\mathcal{B}^*_{(j)}) = y^*/\pi_j,
\]
for every $s_0\ni j$. Whereas, by the zero-invariance property, we have that, for any $s_0 \niton j$,
\[
\lim_{y \to \bm{0}} h_0(s_0)= \lim_{y \to \bm{0}}  h(s_0) =0.
\]

Now, the variance of $e_0$ given the graph $\mathcal{B}^*_{(j)}$ can be decomposed as 
\begin{align*}
\text{V}\left( e_0(s_0; \mathcal{B}^*_{(j)}) \right) &= \sum_{s_0 \niton  j} p(s_0) \left( \theta + h_0(s_0) -\theta-y^*  \right)^2 + \sum_{s_0 \ni j} p(s_0) \left(  e_0(s_0) - \theta-y^*  \right)^2 \\
&=\sum_{s_0 \niton  j} p(s_0) h^2_0(s_0) + (1- \pi_j)(y^*)^2 -2y^*\sum_{s_0 \niton j} p(s_0) h_0(s_0) + \\ 
&\qquad + \sum_{s_0 \ni j} p(s_0) \left(  e_0(s_0) - \theta-y^*  \right)^2.
\end{align*}
Similarly for $e$, such that $\Delta = \text{V}(e_0(s_0; \mathcal{B}^*_{(j)})) - \text{V}(e(s_0; \mathcal{B}^*_{(j)}))$ is given as
\begin{align*}
\Delta &= \sum_{s_0 \ni j} p(s_0) \left(  e_0(s_0) - \theta-y^*  \right)^2 - \sum_{s_0 \ni j} p(s_0) \left(  e(s_0) - \theta-y^*  \right)^2 \\
&\qquad + \sum_{s_0 \niton  j} p(s_0) \left( h_0^2(s_0)-h^2(s_0)\right) -2y^*\sum_{s_0 \niton j} p(s_0) \left( h_0(s_0) -  h(s_0)  \right)
\end{align*}
where, as $y \rightarrow \bm{0}$, the sum of the first two terms above approaches zero because the estimators are elemental, and the third term approaches zero as obtained above.
Thus, for every $j\in F$ and $\epsilon >0$ there exists some $\delta >0 $ such that, for every $\lVert y \rVert < \delta$, the sum of the first three terms are bounded by $\epsilon$ absolutely. We can then write, for $\epsilon^*_{(j)} \in (-\epsilon, \epsilon)$ and $\forall j \in F$,
\[
\text{V}\left( e_0(s_0; \mathcal{B}^*_{(j)}) \right) - \text{V}\left( e(s_0; \mathcal{B}^*_{(j)}) \right)=\epsilon^*_{(j)}-2y^*\sum_{s_0 \niton j} p(s_0) \left( h_0(s_0) -  h(s_0)  \right).
\]

Note that since $e$ differs to $e_0$ given $\mathcal{B}$, there exists some sample $s^*_0$ such that $h_0(s_0^*)\neq h(s_0^*)$ for some value of $y$. Now that $d(s_0)=h(s_0) - h_0(s_0)$ is a non-zero analytic function, all its zeros are isolated, and there exist values of $y$ arbitrarily close to zero where $d(s^*_0) \neq 0.$ Choose a small enough value for $y$ such that $d(s^*_0) \neq 0$ and consider the following.

\textbf{Subcase 2.1:} Suppose that 
$$\sum_{s_0 \niton j} p(s_0) d(s_0) = 0, $$
for all $j \in F.$ Then, as $\sum_{s_0 \in \mathcal{S}_0} p(s_0) d(s_0) = 0 $ due to the unbiasedness, it holds that
\[
\sum_{s_0 \ni j} p(s_0) d(s_0) = 0, \ \forall j \in F.
\]
In other words, $d\in \mathfrak{E}_0\big(\text{Ker}(\tilde{\mathbb{S}}^p) \big).$ By Proposition \ref{prop: ortho-inadmi}, $e_0$ improves $e$ with strict inequality when $d\neq 0$, in which case $\text{V}\left( e \right) \leq \text{V}\left( e_0 \right)$, which yields the contradiction we seek.

\textbf{Subcase 2.2:} Suppose there exists some $j^*\in F$ such that
\[
\sum_{s_0 \niton j^*} p(s_0) d(s_0) \neq 0.
\]
Now, if we take $j^*$ and a small enough value of $y$, we have:
\[
\text{V}\left( e_0(s_0; \mathcal{B}^*_{(j^*)}) \right) - \text{V}\left( e(s_0; \mathcal{B}^*_{(j^*)}) \right)<0 \Leftrightarrow \epsilon^*_{(j^*)} < 2y^*\sum_{s_0 \niton j^*} p(s_0) \left( h_0(s_0) -  h(s_0)  \right).
\]
To find $y^*$ that satisfies the previous inequality, it is enough to let $y^*$ be a sufficiently large positive value if  $\sum_{s_0 \niton j^*} p(s_0) \left( h_0(s_0) -  h(s_0)  \right) > 0$ or a sufficiently large negative value if $\sum_{s_0 \niton j^*} p(s_0) \left( h_0(s_0) -  h(s_0)  \right) < 0$. Thus, there exist some $\overline{y}$ and $\mathcal{B}^*_{(j^*)}$ where $\text{V}\left( e \right) \leq \text{V}\left( e_0\right)$ does not hold, yielding the contraction we seek.
\end{proof}

\noindent 
Proof of Corollary \ref{cor_fullrank}.
\begin{proof}
Note that if $p(s_0)$ is full-rank then by the rank-nullity theorem $\text{Ker}(\tilde{\mathbb{S}}^p)=\text{Ker}(\mathbb{S}^p)=\bf 0.$ By equation (\ref{eq:direct_sum}), $\text{Row}(\mathbb{S}^p)$ is the complete space $\mathbb{R}^{\vert \mathcal{S}_0 \vert}$ and every estimator in $D^{**}$ belongs to $\mathfrak{E}_{\theta}\big( \text{Row}(\mathbb{S}^p) \big)$. It follows then from Theorem \ref{th: main_non_sufficient} that every such estimator is admissible.
\end{proof}

\bibliographystyle{plain}

\end{document}